\documentclass[12pt, reqno]{amsart}
\usepackage{amsmath, amsthm, amscd, amsfonts, amssymb, graphicx}
\usepackage[bookmarksnumbered,plainpages]{hyperref}
\input{mathrsfs.sty}

\newtheorem{theorem}{Theorem}[section]

\newtheorem{corollary}[theorem]{Corollary}
\theoremstyle{definition}

\theoremstyle{remark}
\newtheorem{remark}[theorem]{Remark}

\numberwithin{equation}{section}

\begin{document}

\title{Generalizations of Bohr's inequality in Hilbert $C^*$-modules}
\author[M.S. Moslehian, R. Raji\' c]{Mohammad Sal Moslehian $^1$ and Rajna Raji\' c $^2$}

\address{$^1$ Department of Pure Mathematics, Ferdowsi University of Mashhad,
P.O. Box 1159, Mashhad 91775, Iran.}
\email{moslehian@ferdowsi.um.ac.ir and moslehian@ams.org}

\address{$^2$ Faculty of Mining, Geology and Petroleum
Engineering, University of Zagreb, Pierottijeva 6, 10000 Zagreb,
Croatia} \email{rajna.rajic@zg.t-com.hr}



\begin{abstract}
We present a new operator equality in the framework of Hilbert
$C^*$-modules. As a consequence, we get an extension of the
Euler--Lagrange type identity in the setting of Hilbert bundles as
well as several generalized operator Bohr's inequalities due to
O.~Hirzallah, W.-S. Cheung--J.E.~Pe\v{c}ari\'{c} and F.~Zhang.
\end{abstract}
\maketitle

\section{Introduction}

The classical \emph{Bohr's inequality} states that for any $z, w \in
{\mathbb C}$ and any positive real numbers $r, s$ with
$\frac{1}{r}+\frac{1}{s}=1$,
$$|z+w|^2 \leq r |z|^2 + s |w|^2.$$

\noindent Many interesting generalizations of this inequality have
been obtained in various settings; cf. \cite{C-P, HIR, M-P-F, M-P-P,
P-S, P-D, P-R, RAS1, V-K1, ZHA}. O.~Hirzallah \cite{HIR} showed that
if $A,B$ belong to the algebra ${\mathbb B}({\mathscr H})$ of all
bounded linear operators on a complex Hilbert space ${\mathscr H}$
and $q \geq p>1$ with $1/p + 1/q = 1$, then
$$|A-B|^2 +|(1 - p)A - B|^2 \leq p|A|^2+q|B|^2\,,$$
where $|C|:=(C^*C)^{1/2}$ denotes the absolute value of $C \in
{\mathbb B}({\mathscr H})$.  W.-S.~Cheung and J.E.~Pe\v{c}ari\'{c}
\cite{C-P} extended the above inequality for all positive conjugate
exponents $p, q \in {\mathbb R}$. In addition, F.~Zhang  \cite{ZHA},
among others, generalized the above work of O.~Hirzallah and
presented an identity \cite[Theorem 2]{ZHA} by removing the
condition $q \geq p$.

In this paper we present a new operator equality in the framework of
Hilbert $C^*$-modules. As a consequence, we get a generalization of
the Euler--Lagrange type identity in the setting of Hilbert bundles
over locally compact spaces and extend the operator inequalities of
\cite{C-P}, \cite{HIR} and \cite{ZHA} to get more generalized
inequalities of the Bohr inequality type.


\section{Preliminaries}

The notion of Hilbert $C^*$-module is a generalization of the notion
of Hilbert space. This object was first used by I.~Kaplansky
\cite{KAP}. Hilbert $C^*$-modules are useful tools in Kasparov's
formulation of $KK$-theory, theory of operator algebras, group
representation theory, noncommutative geometry and theory of
operator spaces. It provides a framework for extending the notion of
Morita equivalence to $C^*$-algebras and that of vector bundle to
noncommutative $C^*$-algebras. The theory of Hilbert $C^*$-modules
is interesting on its own right.

Let ${\mathscr A}$ be a $C^*$-algebra and ${\mathscr X}$ be a
complex linear space, which is a right ${\mathscr A}$-module
satisfying $\lambda(xa)=x(\lambda a)=(\lambda x)a$ for $x \in
{\mathscr X},a \in {\mathscr A}, \lambda \in {\mathbb C}$. The
space ${\mathscr X}$ is called a \emph{ (right) pre-Hilbert
$C^*$-module over ${\mathscr A}$} (or a \emph{(right) pre-Hilbert
${\mathscr A}$-module}) if there exists an ${\mathscr A}$-inner
product $\langle .,.\rangle :{\mathscr X} \times {\mathscr X}\to
{\mathscr A}$ satisfying

(i) $\langle x,x\rangle\geq 0$ and $\langle x,x\rangle=0$~~~ if and
only if~~~ $x=0$;

(ii) $\langle x,y+\lambda z\rangle=\langle x,y\rangle+\lambda
\langle x,z\rangle$;

(iii) $\langle x,ya\rangle=\langle x,y\rangle a$;

(iv) $\langle x,y\rangle^*=\langle y,x\rangle$;

\noindent for all $x, y, z \in {\mathscr X},\, \lambda \in {\mathbb
C},\, a \in {\mathscr A}$.

\noindent We can define a norm on ${\mathscr X}$ by $\| x \| =\|
\langle x,x\rangle\| ^\frac{1}{2}$. A pre-Hilbert ${\mathscr
A}$-module is called a \emph{ (right) Hilbert $C^*$-module over
${\mathscr A}$} (or a \emph{(right) Hilbert ${\mathscr A}$-module})
if it is complete with respect to its norm. The notion of \emph{left
Hilbert ${\mathscr A}$-module} can be defined in a similar way.

\noindent Three typical examples of Hilbert $C^*$-modules are as
follows.

(I) Every Hilbert space is a left Hilbert ${\mathbb C}$-module.

(II) Let ${\mathscr A}$ be a $C^*$-algebra. Then ${\mathscr A}$ is a
Hilbert ${\mathscr A}$-module via $\langle a, b\rangle = a^*b \quad
(a, b \in {\mathscr A})$.

(III) Let $$\ell_2({\mathscr A})=\{(a_i)_{i\in{\mathbb N}}:
\sum_{i=1}^\infty a_i^*a_i\textup{ norm-converges in }{\mathscr
A},\,\,a_i\in {\mathscr A}, i=1, 2, \dots\}\,.$$ Then
$\ell_2({\mathscr A})$ is a Hilbert ${\mathscr A}$-module under the
natural operations $\lambda(a_i)+(b_i)=(\lambda a_i+b_i),$ $(a_i)
a=(a_ia)$ and $\langle (a_i),(b_i)\rangle =\sum_{i=1}^\infty
a_i^*b_i$.

\noindent A mapping $T:{\mathscr X}\to {\mathscr Y}$ between Hilbert
${\mathscr A}$-modules is called adjointable if there exists a
mapping $S:{\mathscr Y}\to {\mathscr X}$ such that $\langle
T(x),y\rangle=\langle x,S(y)\rangle$ for all $x\in {\mathscr X},
y\in {\mathscr Y}$. The unique mapping $S$ is denoted by $T^*$ and
is called the adjoint of $T$. It is easy to see that $T$ and $T^*$
must be bounded linear ${\mathscr A}$-module mappings. We denote by
${\mathcal L}({\mathscr X}, {\mathscr Y})$ the space of all
adjointable mappings from ${\mathscr X}$ to ${\mathscr Y}.$ We stand
${\mathcal L}({\mathscr X})$ for the unital $C^*$-algebra ${\mathcal
L}({\mathscr X}, {\mathscr X})$; cf. \cite[p. 8]{LAN}.

\noindent For every $x\in {\mathscr X}$ we define the absolute value
of $x$ as the unique positive square root of $\langle x,x \rangle ,$
that is, $|x|=\langle x,x \rangle ^\frac{1}{2}$. We refer the reader
to \cite{MUR} for undefined notions on $C^*$-algebra theory and to
\cite{LAN, M-T, R-W} for more information on Hilbert $C^*$-modules.

Throughout the paper, we assume that ${\mathscr X}$ and ${\mathscr
Y}$ are Hilbert $C^*$-modules over a unital $C^*$-algebra ${\mathscr
A}$ with unit $e$. The identity operator on a set $E$ is denoted by
$I_E$. By ${\mathcal Z}({\mathscr A})$ we denote the center of a
$C^*$-algebra ${\mathscr A}.$


\section{Main results}

\begin{theorem}\label{prvi}
Let $\alpha, \beta, \gamma \in {\mathbb R}$ and let $T,S\in
{\mathcal L}({\mathscr X}, {\mathscr Y})$ be such that $T^*S \in
{\mathcal L}({\mathscr X})$ is self-adjoint and $\alpha T^*T+\beta
S^*S=\gamma I_{\mathscr X}$. Then
$$\alpha\beta |Tx+Sy|^2+|\beta Sx-\alpha Ty|^2 = \beta\gamma |x|^2+\alpha\gamma |y|^2 $$
for all $x,y\in {\mathscr X}$.
\end{theorem}
\begin{proof}
Since
$$\begin{array}{rcl}
|Tx+Sy|^2 & = & \langle Tx+Sy, Tx+Sy \rangle \\
& = & \langle Tx,Tx \rangle +\langle Tx,Sy \rangle +\langle Sy,Tx
\rangle +\langle Sy,Sy \rangle \\
& = & \langle T^*Tx,x \rangle +\langle S^*Tx,y \rangle +\langle
T^*Sy,x \rangle +\langle S^*Sy,y \rangle ,\\
\end{array}$$
and
$$\begin{array}{rcl}
|\beta Sx-\alpha Ty|^2 & = & \langle \beta Sx-\alpha Ty, \beta
Sx-\alpha Ty \rangle \\
& = & \beta^2\langle Sx,Sx \rangle -\alpha\beta \langle Sx,Ty
\rangle -\alpha\beta \langle Ty,Sx \rangle +\alpha^2 \langle Ty,Ty
\rangle \\
& = & \beta^2 \langle S^*Sx,x \rangle -\alpha\beta \langle
T^*Sx,y\rangle -\alpha\beta \langle S^*Ty,x \rangle +\alpha^2
\langle T^*Ty,y \rangle \\
\end{array}$$
we get
$$\begin{array}{rcl}
\alpha\beta |Tx+Sy|^2+|\beta Sx-\alpha Ty|^2 & = & \beta \langle
(\alpha T^*T+\beta S^*S)x,x \rangle +\alpha \langle (\alpha
T^*T+\beta S^*S)y,y \rangle \\
& = & \beta\gamma |x|^2+\alpha\gamma |y|^2\,. \\
\end{array}$$
\end{proof}


The following useful result is deduced from Theorem~\ref{prvi}. It
can be proved directly as well.

\begin{theorem}\label{cprvi}
Let $\alpha, \beta, \gamma \in {\mathbb R}$ and let $x, y\in
{\mathscr X}$ be such that $\langle x, y\rangle$ is self-adjoint and
$\alpha \langle x, x\rangle + \beta \langle y, y\rangle=\gamma e$.
Then
$$\alpha\beta |xa+yb|^2+|\beta ya-\alpha
xb|^2=\beta\gamma|a|^2+\alpha\gamma|b|^2$$ for all $a,b\in {\mathscr
A}$.
\end{theorem}
\begin{proof}
For each $z\in {\mathscr X},$ consider two mappings $T_z: {\mathscr
A} \to {\mathscr X}$ and $\hat{z}: {\mathscr X} \to {\mathscr A}$
defined by $T_z(a)=za$ and $\hat{z}(v)= \langle z, v\rangle$. The
adjoint of $T_z$ is $\hat{z}$, since
$$\langle T_z(a), v\rangle=\langle za, v\rangle=a^*\langle z, v\rangle=a^*\hat{z}(v)$$
for all $a\in {\mathscr A}$ and $v\in {\mathscr X}.$ Furthermore,
$T_x^*T_y$ is self-adjoint since
$$T_x^*T_y(a)=\hat{x}(ya)=\langle x,ya\rangle =\langle x,y\rangle
a=\langle y,x\rangle a=\langle y,xa\rangle
=\hat{y}(xa)=T_y^*T_x(a)$$ for all $a\in {\mathscr A}.$ In addition,
\begin{eqnarray*}
\left(\alpha T_x^*T_x+\beta T_y^*T_y\right)(a)&=&
\alpha\hat{x}(xa)+\beta\hat{y}(ya)\\
&=&\alpha \langle x, xa\rangle+ \beta\langle y, ya\rangle\\
&=& (\alpha\langle x, x\rangle+\beta\langle y, y\rangle)a\\
&=& \gamma I_{\mathscr A}(a)
\end{eqnarray*}
for all $a\in {\mathscr A}.$ Applying Theorem \ref{prvi} with
$T=T_x$ and $S=T_y$, we obtain
$$\alpha\beta |T_xa+T_yb|^2 + |\beta T_ya-\alpha
T_xb|^2=\beta\gamma |a|^2+ \alpha\gamma |b|^2\qquad (a, b\in
{\mathscr A})$$ which proves the theorem.
\end{proof}


Applying Theorem~\ref{cprvi} for elements of the Hilbert
$C^*$-module $\ell_2({\mathscr A})$ we get the following result.

\begin{corollary}
Let $\alpha, \beta, \gamma \in {\mathbb R}$ and let $(a_i)_i,
(b_i)_i \in \ell_2({\mathscr A})$ be such that $\sum_{i=1}^\infty
a_i^*b_i$ is self-adjoint and $\alpha \sum_{i=1}^\infty
|a_i|^2+\beta \sum_{i=1}^\infty |b_i|^2=\gamma e.$ Then
$$\alpha\beta \sum_{i=1}^\infty |a_ia+b_ib|^2 + \sum_{i=1}^\infty |\beta b_ia-\alpha
a_ib|^2 = \beta\gamma |a|^2 + \alpha\gamma |b|^2$$ for all $a,b\in
{\mathscr A}$.
\end{corollary}


Recall that the space ${\mathbb B}(\mathscr H,\mathscr K)$ of all
bounded linear operators between Hilbert spaces $\mathscr H$ and
$\mathscr K$ can be regarded as a Hilbert $C^*$-module over the
$C^*$-algebra ${\mathbb B}(\mathscr H)$ via $\langle
T,S\rangle=T^*S$. Then the direct sum
$${\mathscr X}=\underbrace{{\mathbb B}(\mathscr H,\mathscr K)\oplus \cdots \oplus {\mathbb B}(\mathscr H,\mathscr K)}_{n}
=\{(T_1,\dots,T_n): T_i\in {\mathbb B}(\mathscr H,\mathscr K),
i=1,\dots,n\}$$ is a Hilbert ${\mathbb B}(\mathscr H)$-module, where
the inner product is defined as $\langle (T_i)_i, (S_i)_i \rangle
=\sum_{i=1}^nT_i^*S_i,$ $(T_i, S_i\in {\mathbb B}(\mathscr
H,\mathscr K))$. For such a Hilbert $C^*$-module ${\mathscr X}$,
Theorem \ref{cprvi} can be stated as follows.

\begin{corollary}
Let $\alpha, \beta, \gamma \in {\mathbb R}$ and let
$T_1,\dots,T_n,S_1,\dots,S_n\in {\mathbb B}(\mathscr H,\mathscr K)$
be such that $\sum_{i=1}^nT_i^*S_i$ is self-adjoint and $\alpha
\sum_{i=1}^n T_i^*T_i +\beta \sum_{i=1}^nS_i^*S_i =\gamma
I_{\mathscr H}$. Then
$$\alpha\beta \sum_{i=1}^n(T_iA+S_iB)^*(T_iA+S_iB) + \sum_{i=1}^n(\beta
S_iA-\alpha T_iB)^*(\beta S_iA-\alpha T_iB) = \beta\gamma |A|^2 +
\alpha\gamma |B|^2$$ for all $A,B\in {\mathbb B}(\mathscr H)$.
\end{corollary}


As another consequence of Theorem \ref{prvi} we have the following
generalization of the Euler--Lagrange type identity \cite{T-R-S-T}
in the framework of Hilbert $C^*$-modules.

\begin{theorem}\label{eul-lagr}
Let $\alpha, \beta, \gamma \in {\mathbb R}$ and let $a,b\in
{\mathcal Z}({\mathscr A})$ be such that $a^*b$ is self-adjoint and
$\alpha a^*a+\beta b^*b=\gamma e$. Then
$$\alpha\beta |xa+yb|^2 + |\beta xb-\alpha
ya|^2 = \beta\gamma |x|^2 + \alpha\gamma |y|^2$$ for all $x,y \in
{\mathscr X}$.
\end{theorem}
\begin{proof}
For each $c\in {\mathcal Z}({\mathscr A})$, the mapping
$T_c:{\mathscr X}\rightarrow {\mathscr X}$ defined by $T_c(x)=xc$
has the adjoint $T_{c^*}$, since
$$\langle T_c(x),y \rangle =\langle xc,y \rangle =c^*\langle x,y
\rangle =\langle x,y \rangle c^* =\langle x,yc^* \rangle =\langle
x,T_{c^*}(y) \rangle \qquad (x, y \in {\mathscr X})\,.$$
Furthermore,
$$T_a^*T_b=T_{a^*}T_b=T_{a^*b}=T_{b^*a}=T_{b^*}T_a=T_b^*T_a,$$
and
$$\alpha T_a^*T_a+\beta T_b^*T_b=T_{\alpha a^*a+\beta
b^*b}=T_{\gamma e}=\gamma T_e=\gamma I_{\mathscr X}.$$ The result
follows by applying Theorem~\ref{prvi} with $T=T_a$ and $S=T_b$.
\end{proof}

Next consider a locally compact space $K$ and assume that ${\mathcal
B}=\bigcup_{t \in K}{\mathscr H}_t$ is a bundle of Hilbert spaces
over $K$ which satisfies appropriate continuous properties. Then the
set $C_0(K,{\mathcal B})$ of all continuous mappings $\varphi: K \to
{\mathcal B}$ which vanishes at infinity and fulfills $\varphi(t)
\in {\mathscr H}_t\,\,(t \in K)$ is a Hilbert $C_0(K)$-module via
$$(\varphi \cdot f)(t):=\varphi(t)f(t) \qquad
\langle \varphi, \psi\rangle(t):=\langle \psi(t)|\varphi(t)\rangle_t
\qquad (f \in C_0(K), \varphi, \psi \in C_0(K,{\mathcal B}))\,,$$
where $\langle \cdot |\cdot\rangle_t$ stands for the inner product
of Hilbert space ${\mathscr H}_t$ (we denote the induced norm in
${\mathscr H}_t$ by $\|.\|_t$); see \cite{KAP}. It is easy to see
that $|\varphi|$ is the function $t \mapsto \|\varphi(t)\|_t\,\,(t
\in K)$. The following result immediately follows from
Theorem~\ref{eul-lagr}. The special case, where $K$ is singleton,
gives rise to the classical Euler--Lagrange type identity.

\begin{corollary}[Generalized Euler--Lagrange type identity] Suppose
that ${\mathcal B}=\bigcup_{t \in K}{\mathscr H}_t$ is a bundle of
Hilbert spaces over a locally compact space $K$. Let $\alpha, \beta,
\gamma \in {\mathbb R}$ and let $f, g\in C_0(K)$ be real functions
such that $\alpha f(t)^2 + \beta g(t)^2=\gamma\,\,(t \in K)$. Then
\begin{eqnarray*}
\sup_{t\in K} \big(\alpha\beta \|\varphi(t)f(t)+ \psi(t)g(t)\|_t^2 +
\|\beta \varphi(t)g(t) - \alpha \psi(t)f(t)\|_t^2\big) = \sup_{t \in
K}\big( \beta\gamma \|\varphi(t)\|_t^2 + \alpha\gamma
\|\psi(t)\|_t^2\big)
\end{eqnarray*}
for all $\varphi, \psi \in C_0(K,{\mathcal B})$.
\end{corollary}


A generalization of Bohr's inequality in Hilbert $C^*$-modules
deduced from Theorem~\ref{eul-lagr} as well as some generalizations
of main results of \cite[Theorem 1 and Corollary 1]{C-P} and
\cite[Theorem 1]{HIR} can be presented in the following theorem.

\begin{theorem}\label{3.7}
Let $p,q>1$ be conjugate components. Then
\begin{eqnarray}\label{xyp}
|x-y|^2+\frac{1}{p-1}|(1-p)x-y|^2=p|x|^2+q|y|^2
\end{eqnarray}
for all $x,y$ in a Hilbert $C^*$-module ${\mathscr X}$. Moreover,
\begin{eqnarray*}
&&(i)\qquad\,\,  |x-y|^2 +|(1 - p)x - y|^2 \leq p|x|^2+q|y|^2 \Leftrightarrow p \leq 2 {\rm ~or~ } (1- p)x=y\,,\\
&&(ii)\qquad\,\, |x-y|^2 +|(1 - p)x - y|^2 \geq p|x|^2+q|y|^2 \Leftrightarrow p \geq 2 {\rm ~or~ } (1- p)x=y\,.\\
\end{eqnarray*}
Furthermore, in (i) and (ii) the equality holds on the left hand
side of the equivalence if and only if $p=q=2$ or $(1-p)x=y$.
\end{theorem}
\begin{proof}
Apply Theorem~\ref{eul-lagr} with $a=e,$ $b=-e,$ $\alpha =1,$ $\beta
=p-1$ and $\gamma =p$.

\noindent To achieve (i), use \eqref{xyp} and the fact that $1 \leq
\frac{1}{p-1}$ if and only if $p \leq 2$. The rest can be proved in
a similar way.
\end{proof}
\begin{remark}
A similar assertion can be proved for the case where $p<1$, see
\cite[Theorem 2]{C-P}. Theorem \ref{3.7} generalizes Bohr's
inequality in Hilbert $C^*$-modules only when $p$ and $q$ are
positive conjugate exponents and $p\leq 2$. Interchanging $x
\leftrightarrow y$ and $p \leftrightarrow q$ in \eqref{xyp} we also
have
$$|x-y|^2+\frac{1}{q-1} |(1-q)y-x|^2 =p|x|^2+g|y|^2.$$
From this we have a generalization of Bohr's inequality in the case
when $p>2$ (that is, $1<q<2$). Namely, statement (i) of Theorem
\ref{3.7} now reads as follows:
$$|x-y|^2+|(1-q)y-x|^2\le p|x|^2+q|y|^2 \Leftrightarrow p \geq 2 {\rm ~or~ } (1-q)y=x\,.$$
\end{remark}

One more consequence of Theorem~\ref{prvi} is the following result.

\begin{corollary}\label{bohr2}
Let $\alpha, \beta$ be positive real numbers satisfying $\alpha
+\beta =1$. Let $T,S\in {\mathcal L}({\mathscr X}, {\mathscr Y})$ be
such that $T^*S\in {\mathcal L}({\mathscr X})$ is self-adjoint and
$\alpha T^*T+\beta S^*S=I_{\mathscr X}$. Then
$$|\beta Sx+\alpha Ty|^2\le \beta |x|^2+\alpha |y|^2$$
for all $x,y\in \mathscr X$.
\end{corollary}

\begin{proof}
The result follows immediately from Theorem~\ref{prvi} by taking
$\gamma =\alpha +\beta =1$.
\end{proof}

Our next result is a generalization of Corollary~\ref{bohr2} in the
case of an arbitrary number of finitely many elements of ${\mathcal
L}({\mathscr X})$.

\begin{theorem}\label{bohrn}
Let $n\ge 2$ be a positive integer, let $T_1,\dots,T_n\in {\mathcal
L}({\mathscr X})$, let $T_1^*T_2$ be self-adjoint, and let
$t_1,\dots,t_n$ be positive real numbers such that $\sum_{i=1}^n
t_i=1$ and $\sum_{i=1}^nt_i|T_i|^2=I_{\mathscr X}$. For $n \geq 3$,
assume $T_1$ or $T_2$ is invertible in ${\mathcal L}({\mathscr X})$,
operators $T_3,\dots,T_n$ are self-adjoint, and $T_i|T_j|=|T_j|T_i$
for all $1 \leq i<j \leq n$. Then
\begin{eqnarray}\label{zhang}
|t_1T_1x_1+\cdots +t_nT_nx_n|^2\le t_1|x_1|^2+\cdots +t_n|x_n|^2
\end{eqnarray}
for all $x_1,\dots,x_n\in \mathscr X$.
\end{theorem}

\begin{proof}
We prove the statement by induction on $n$. The base case is true by
Corollary~\ref{bohr2}.

\noindent Suppose now that the inequality \eqref{zhang} holds for
$n-1\ge 2$ elements and let us show that it holds for $n$ elements
$T_1, \dots, T_n$. Assume that $T_k$ is invertible for some
$k\in\{1,2\}$. Then
$$I_{\mathscr X}-t_n|T_n|^2=\sum_{i=1}^{n-1}t_i|T_i|^2 \ge
t_k|T_k|^2 > 0$$ implies invertibility of $I_{\mathscr
X}-t_n|T_n|^2$ in ${\mathcal L}({\mathscr X})$. Let us put
$$y=\sum_{i=1}^{n-1}s_iS_ix_i,$$
where
$$s_i=\frac{t_i}{1-t_n}, \quad S_i=\sqrt{1-t_n} T_i(I_{\mathscr X}-t_n|T_n|^2)^{-\frac{1}{2}} \quad
(i=1,\dots, n-1).$$ Clearly, $\sum_{i=1}^{n-1}s_i=1,$ $S_i$ are
self-adjoint for $i=3,\dots,n-1,$ $S_k$ is invertible as the product
of two invertible operators, and $S_1^*S_2$ is self-adjoint since
$$\begin{array}{rcl}
S_1^*S_2 & = &
(1-t_n)(I_{\mathscr X}-t_n|T_n|^2)^{-\frac{1}{2}}T_1^*T_2(I_{\mathscr X}-t_n|T_n|^2)^{-\frac{1}{2}}\\
& = & (1-t_n)(I_{\mathscr
X}-t_n|T_n|^2)^{-\frac{1}{2}}T_2^*T_1(I_{\mathscr
X}-t_n|T_n|^2)^{-\frac{1}{2}} = S_2^*S_1.\\
\end{array}$$
Observe that for each $1 \leq i<j \leq n$, $T_i|T_j|=|T_j|T_i$
implies $|T_i||T_j|=|T_j||T_i|$ from which we get
$$\begin{array}{rcl}
|S_j|^2 = (1-t_n)(I_{\mathscr
X}-t_n|T_n|^2)^{-\frac{1}{2}}|T_j|^2(I_{\mathscr
X}-t_n|T_n|^2)^{-\frac{1}{2}} =  (1-t_n)|T_j|^2(I_{\mathscr
X}-t_n|T_n|^2)^{-1}
\end{array}$$
for all $1 \leq j \leq n-1$. It follows that
$$\begin{array}{rcl}
S_i|S_j|^2 & = & (1-t_n)^\frac{3}{2}T_i(I_{\mathscr
X}-t_n|T_n|^2)^{-\frac{1}{2}}|T_j|^2(I_{\mathscr
X}-t_n|T_n|^2)^{-1}\\
& = & (1-t_n)^\frac{3}{2}|T_j|^2(I_{\mathscr
X}-t_n|T_n|^2)^{-1}T_i(I_{\mathscr X}-t_n|T_n|^2)^{-\frac{1}{2}} = |S_j|^2S_i\,,\\
\end{array}$$
whence $S_i|S_j|=|S_j|S_i$ for all $1 \leq i<j \leq n-1$. In
addition,
$$\begin{array}{rcl}
\displaystyle\sum_{i=1}^{n-1}s_i|S_i|^2 =
\sum_{i=1}^{n-1}t_i|T_i|^2(I_{\mathscr X}-t_n|T_n|^2)^{-1} =
(I_{\mathscr X}-t_n|T_n|^2)(I_{\mathscr X}-t_n|T_n|^2)^{-1} =
I_{\mathscr X}.
\end{array}$$
By the inductive assumption we conclude that
\begin{equation}\label{y}
|y|^2\le \sum_{i=1}^{n-1}s_i|x_i|^2.
\end{equation}
Let us denote $W=\frac{1}{\sqrt{1-t_n}}(I_{\mathscr
X}-t_n|T_n|^2)^\frac{1}{2}.$ Then
\begin{equation}\label{uv}
(1-t_n)|W|^2+t_n|T_n|^2=(1-t_n)\bigg|\frac{1}{\sqrt{1-t_n}}
(I_{\mathscr
X}-t_n|T_n|^2)^\frac{1}{2}\bigg|^2+t_n|T_n|^2=I_{\mathscr X}.
\end{equation}
Also,
\begin{eqnarray}\label{pom}
(1-t_n)Wy & = & \displaystyle\sqrt{1-t_n}(I_{\mathscr
X}-t_n|T_n|^2)^\frac{1}{2}\sum_{i=1}^{n-1}s_iS_ix_i \nonumber\\
& = & \displaystyle\sqrt{1-t_n}(I_{\mathscr
X}-t_n|T_n|^2)^\frac{1}{2}\sum_{i=1}^{n-1}s_i\sqrt{1-t_n}T_i(I_{\mathscr
X}-t_n|T_n|^2)^{-\frac{1}{2}}x_i \nonumber\\
& = & \displaystyle\sum_{i=1}^{n-1}t_iT_ix_i.
\end{eqnarray}
Using (\ref{y}), (\ref{uv}), (\ref{pom}) and the fact that $T_n$ is
self-adjoint, we obtain
$$\begin{array}{rcl}
\displaystyle|t_1T_1x_1+\cdots +t_{n-1}T_{n-1}x_{n-1}+t_nT_nx_n|^2 &
= &
\displaystyle|(1-t_n)Wy+t_nT_nx_n|^2\\
& \le & (1-t_n)|y|^2+t_n|x_n|^2\\
& \le & t_1|x_1|^2+\cdots +t_{n-1}|x_{n-1}|^2+t_n|x_n|^2.\\
\end{array}$$
This proves the theorem.
\end{proof}

\begin{remark}
In Theorems \ref{prvi} and \ref{bohrn} and Corollary \ref{bohr2} we
do not need the underlying $C^*$-algebra to be unital.
\end{remark}

\begin{corollary}\label{bohrncor}
Let $n\ge 2$ be a positive integer, let $a_1,\dots,a_n\in {\mathcal
Z}({\mathscr A})$, let $a_1^*a_2$ be self-adjoint, and let
$t_1,\dots,t_n$ be positive real numbers such that $\sum_{i=1}^n
t_i=1$ and $\sum_{i=1}^nt_i|a_i|^2=e$. For $n\geq 3$, assume that
$a_1$ or $a_2$ is invertible in ${\mathscr A}$ and $a_3,\dots,a_n$
are self-adjoint. Then
\begin{eqnarray}\label{zhangmain}
|t_1x_1a_1+\cdots +t_nx_na_n|^2\le t_1|x_1|^2+\cdots +t_n|x_n|^2
\end{eqnarray}
for all $x_1,\dots,x_n\in \mathscr X$.
\end{corollary}

\begin{proof}
Let $T_i:{\mathscr X}\rightarrow {\mathscr X}$ be the mapping
defined by $T_i(x)=xa_i$ ($i=1,\dots,n$). Then, as seen in the proof
of Theorem \ref{eul-lagr}, $T_i^*=T_{a_i^*}$ and so
$|T_i|=T_{|a_i|}$. One can easily verify that the operators $T_i$
satisfy the assumptions of Theorem~\ref{bohrn}. Thus inequality
(\ref{zhang}) turns into (\ref{zhangmain}).
\end{proof}

When ${\mathscr X}$ is the $C^*$-algebra ${\mathbb B}({\mathscr H})$
regarded as a Hilbert $C^*$-module over itself, and $a_i=I_{\mathscr
H}\,\, (i=1, \dots, n)$, the above result reduces to Theorem 7 of
Zhang \cite{ZHA}, which is an AM-QM operator inequality.

\begin{corollary}
Let $A_1, \dots, A_n \in {\mathbb B}({\mathscr H})$. Then, for any
set of nonnegative numbers $t_1, \dots, t_n$ with
$\sum_{i=1}^nt_i=1$,
\begin{eqnarray*}
\left|\sum_{i=1}^nt_iA_i\right|^2 \leq \sum_{i=1}^nt_i|A_i|^2 \,.
\end{eqnarray*}
\end{corollary}

\textbf{Acknowledgement.} This research was in part supported by a
grant from Center for Research in Modeling and Computation of Linear
and Nonlinear Systems (CRMCS).


\end{document}